\documentclass[11pt]{amsart}
\usepackage{amsfonts}
\usepackage{amsmath,amsthm,latexsym,amssymb,graphicx}
\usepackage{graphicx,color,enumerate}

\numberwithin{equation}{section}

\newtheorem{theorem}{Theorem}

\newtheorem{remark}{Remark}

\numberwithin{theorem}{section}
\numberwithin{corollary}{section}
\numberwithin{lemma}{section}
\numberwithin{definition}{section}
\numberwithin{proposition}{section}
\numberwithin{remark}{section}

\newcommand{\R}{\mathbb R}

\newcommand{\medint}{-\kern  -,375cm\int}

\def\deg{\mathrm{deg}}
\thanks{
$^{*}$ Dipartimento di Matematica e Applicazioni ``R. Caccioppoli'', Universit\`{a}
degli Studi di Napoli ``Federico II'', Complesso Monte S. Angelo, via Cintia
- 80126 Napoli, Italy; email: brandolini@unina.it; fchiacch@unina.it
}

\begin{document}
\title{A  remark on the radial minimizer of the Ginzburg-Landau functional}
\author{Barbara Brandolini$^{*}$}
\author{Francesco Chiacchio$^{*}$}

\maketitle

{\Small \textsc{Abstract}. 
Denote by $E_\epsilon$ the Ginzburg-Landau functional in the plane and let $\tilde  u_\varepsilon$ be the radial solution to the Euler equation associated to the problem $\min \left\{E_\varepsilon(u,B_1): \>\left.  u\right\vert _{\partial B_{1}}=(\cos \vartheta,\sin \vartheta)\right\}$. Let $\Omega\subset \R^2$ be a smooth, bounded domain with the same area as $B_1$. Denoted by 
$$\mathcal{K}=\left\{v=(v_1,v_2) \in H^1(\Omega;\R^2):\> \int_\Omega v_1\,dx=\int_\Omega v_2\,dx=0,\> \int_\Omega |v|^2\,dx\ge \int_{B_1} |\tilde u_\varepsilon|^2\,dx\right\},$$
 we prove 
$$
\min_{v \in \mathcal{K}} E_\varepsilon (v,\Omega)\le E_\varepsilon (\tilde u_\varepsilon,B_1).
$$
}


\section{Introduction}

The Ginzburg-Landau energy has as order parameter a vectorial field $u\in H^1(\Omega;\R^2)$  and is defined as
\begin{equation*}\label{E}
E_\varepsilon(u,\Omega)=\frac{1}{2}\int_\Omega |\nabla u|^2\,dx+\frac{1}{4\varepsilon^2}\int_\Omega \left(|u|^2-1\right)^2\,dx,
\end{equation*}
where $\Omega\subset\R^2$ is a bounded domain and  $\varepsilon>0$.  
This kind of functionals has been originally introduced as a phenomenological phase-field type free-energy of a superconductor, near the superconducting transition, in absence of an external magnetic field. 
Moreover these functionals have been used in superfluids such as Helium II. In this context $u$ represents the wave function of the superflluid part of liquid and the parameter  $\varepsilon$, which has the dimension of a length, depends on the material and its temperature.  
The Ginzburg-Landau functionals have deserved a great attention by the mathematical community too.  Starting from the paper \cite{BBH} by  Bethuel, Brezis and H\'elein, many mathematicians have been interested in studying minimization problems for the Ginzburg-Landau energy with several constraints, also because, besides the physical motivation, these problems appear as the simplest nontrivial examples of vector field minimization problems.

In \cite{BBH} the authors consider Dirichlet boundary conditions $g\in C^1(\partial\Omega; \mathbb{S}^1)$ and study the asymptotic behavior, as $\varepsilon \to 0$,  of minimizers $u_\varepsilon$, which clearly satisfies the following Euler equation 
\begin{equation}\label{pb}
\left\{\begin{array}{ll}
-\Delta u_\varepsilon =\frac{1}{\varepsilon^2}u_\varepsilon\left(1-|u_\varepsilon|^2\right) & \mbox{in}\> \Omega  \\\\ u_\varepsilon=g & \mbox{on}\> \partial\Omega.
\end{array}
\right.
\end{equation}
It turns out that the value $d = \deg(g,\partial\Omega)$ (i.e., the Brouwer degree or winding number of $g$ considered as a map from $\partial\Omega$ into $\mathbb{S}^1$) plays a crucial role in the asymptotic analysis of $u_\varepsilon$. 

In the case $\Omega=B_1$ (the unit ball in $\R^2$ centered at the origin), $g(x)=x$, it is natural to look for radial solutions to \eqref{pb}. Indeed, in \cite{G,BBHbook,HH} the authors prove, among other things, that 
 \eqref{pb} has a unique radial solution, that is a solution of the form
\begin{equation}\label{u}
\tilde u_\varepsilon(x)=\tilde f_\varepsilon(|x|)\left(\cos \vartheta, \sin \vartheta\right)
\end{equation}
with $\tilde f_\varepsilon \ge 0$. Moreover $\tilde f_\varepsilon'>0$; thus, summarizing, $\tilde f_\varepsilon$ is the solution to the following problem 
\begin{equation}\label{f}
\left\{\begin{array}{ll}
-\tilde f_\varepsilon''-\dfrac{\tilde f_\varepsilon'}{r}+\dfrac{\tilde f_\varepsilon}{r^2}=\dfrac{1}{\varepsilon^2}\tilde f_\varepsilon\left(1-\tilde f_\varepsilon^2\right) & \mbox{in}\> (0,1)
\\ \\
\tilde f_\varepsilon(0)=0,\> \tilde f_\varepsilon(1)=1,\> \tilde f_\varepsilon\ge 0,\> \tilde f_\varepsilon'>0.
\end{array}
\right.
\end{equation} 
It is conjectured that the radial solution \eqref{u} is the unique minimizer of $E_\varepsilon$ on $B_1$.
In \cite{M} (see also \cite{LL}) the author  gives a partial answer to such a conjecture, proving that $\tilde u_\varepsilon$ is stable, in the sense that the quadratic form associated to $E_\varepsilon(\tilde u_\varepsilon,B_1)$ is positive definite. 

Other types of boundary conditions, for instance prescribed degree boundary conditions,  have been considered in \cite{BR, D}. 

In this paper we let $\Omega$ vary among smooth domains with fixed area and prove that  
 the map $\tilde u_\varepsilon$ in \eqref{u} provides an upper bound for the energy $E_\varepsilon$ on the class $\mathcal{K}$ we are going to introduce. 
 
 \begin{theorem}\label{1}
Let $\varepsilon>0$ and $\Omega\subset \R^2$ be a smooth, bounded domain such that $|\Omega|=|B_1|$. Denoted by
$$\mathcal{K}=\left\{v=(v_1,v_2) \in H^1(\Omega;\R^2):\> \int_\Omega v_1\,dx=\int_\Omega v_2\,dx=0,\> \int_\Omega |v|^2\,dx\ge \int_{B_1} |\tilde u_\varepsilon|^2\,dx\right\},$$
it holds
\begin{equation}\label{max}
 \min_{v \in \mathcal{K}} E_\varepsilon (v,\Omega)\le E_\varepsilon (\tilde u_\varepsilon,B_1).
\end{equation} 
\end{theorem}

\section{Proof of Theorem \ref{1}}
 
Define  the following continuous extension of $\tilde f_\varepsilon$
$$
 f_\varepsilon(r)=\left\{\begin{array}{ll}
\tilde f_\varepsilon(r) & \mbox{if} \>\> 0\le r \le 1 \\ \\ 1 &\mbox{if}\>\>  r>1
\end{array}
\right.
$$
and the correspondent   vector field extending $\tilde u_\varepsilon$ to the whole $\R^2$
$$\phi_\varepsilon(x)=\left(\phi_{\varepsilon,1}(x),\phi_{\varepsilon,2}(x)\right)=f_\varepsilon(|x|)\left(\cos\vartheta, \sin\vartheta\right).$$
It is possible (see  \cite{W}, see also \cite{AB}) to choose the origin in such a way that 
\begin{equation}\label{orth}
\int_\Omega \phi_{\varepsilon,1}\,dx=\int_\Omega \phi_{\varepsilon,2}\,dx=0.
\end{equation}
Note that $\phi_\varepsilon\in \mathcal{K}$. Indeed, besides \eqref{orth}, it holds
$$
\int_\Omega |\phi_\varepsilon|^2\,dx = \int_{\Omega \cap B_1} |\tilde u_\varepsilon|^2\,dx+|\Omega \setminus B_1| \ge \int_{B_1} |\tilde u_\varepsilon|^2\,dx,
$$
since $|\tilde u_\varepsilon| \le 1$ in $B_1$. A direct computation yields
\begin{eqnarray*}
E_\varepsilon(\phi_\varepsilon,\Omega)&=&\frac{1}{2}\int_\Omega \left(f_\varepsilon'(|x|)^2+\frac{ f_\varepsilon(|x|)^2}{|x|^2}\right)\,dx+\frac{1}{4\varepsilon^2}\int_\Omega \left(f_\varepsilon(|x|)^2-1\right)^2\,dx
\\
&=&\int_\Omega B_\varepsilon(|x|)\,dx,
\end{eqnarray*}
where 
$$
B_\varepsilon(r)=\frac{1}{2}\left(f_\varepsilon'(r)^2+\frac{f_\varepsilon(r)^2}{r^2}\right)+\frac{1}{4\varepsilon^2}\left(f_\varepsilon(r)^2-1\right)^2.
$$
Using \eqref{f} it is straightforward to verify that  
\begin{equation*}\label{decreasing}
B_\varepsilon'(r)=-\frac{2}{\varepsilon^2}f_\varepsilon(r)f_\varepsilon'(r)\left(1-f_\varepsilon(r)^2\right)-\frac{1}{r}\left(f_\varepsilon'(r)-\frac{f_\varepsilon(r)}{r}\right)^2, \quad 0<r<1,
\end{equation*}
while, when $r>1$, it holds $B_\varepsilon(r)=\frac{1}{2r^2}$. Thus $B_\varepsilon(r)$ is a decreasing function in $(0,+\infty)$. By Hardy-Littlewood inequality (see for instance \cite{HLP}) we finally get
$$
E_\varepsilon(\phi_\varepsilon,\Omega)=\int_\Omega B_\varepsilon(|x|) \,dx\le \int_{B_1} B_\varepsilon(|x|)\,dx=E_\varepsilon(\tilde u_\varepsilon,B_1)
$$
and hence \eqref{max}.

\begin{remark}
The appearance of the function $\tilde u_\varepsilon$ (i.e., the candidate to be the unique minimizer of $E_\varepsilon$ in $B_1$ under the Dirichlet boundary condition $g(x)=x$) in \eqref{max} as an  upper bound  of the energy in the class $\mathcal{K}$ could seem odd. On the other hand such a phenomenon  looks more natural once one notices an analogy between the problem under consideration and the maximization problem of the first nontrivial eigenvalue $\mu_1(\Omega)$ of the Neumann Laplacian among sets with prescribed area. 
As well-known, if $\Omega$ is a smooth, bounded domain of $\R^2$, $\mu_1(\Omega)$ can be variationally characterized as
$$
\mu_1(\Omega)=\left\{\int_\Omega |\nabla z|^2 : z \in H^1(\Omega;\R),\> \int_\Omega |z|^2\,dx=1,\> \int_\Omega z \,dx=0\right\}.
$$
If $|\Omega|=|B_1|$ the celebrated Szeg\"o-Weinberger inequality in the plane  (see \cite{W}, see also \cite{S,B,AB,H,GP,CdB}) states 
\begin{equation}\label{sw}
\mu_1(\Omega) \le \mu_1(B_1).
\end{equation}
Moreover, $\mu_1(B_1)$ is achieved by the functions $J_1(j_{1,1}'|x|)\cos \vartheta$ or $J_1(j_{1,1}'|x|)\sin \vartheta$, where $J_1$ is the Bessel function of the first kind and $j_{1,1}'$ is the first zero of its derivative. 
 The role played by $J_1$ in \eqref{sw} is now played by the function $\tilde f_\varepsilon$. 

\end{remark}

\end{document}